 \documentclass[letterpaper, paper,11pt]{AAS}		

\usepackage{multirow}
\usepackage{bm}
\usepackage{amsmath}
\usepackage{amssymb}
\usepackage{physics}

\usepackage{hyperref}
\hypersetup{colorlinks,allcolors=black}
\usepackage{overcite}
\usepackage{footnpag}			      	
\usepackage{siunitx}
\usepackage{float}
\usepackage{inputenc}
\usepackage{leftidx}
\usepackage{booktabs}
\usepackage{romannum}

\usepackage{algorithm,algcompatible,algorithmicx}
\algnewcommand\INPUT{\item[\textbf{Input:}]}%
\algnewcommand\OUTPUT{\item[\textbf{Output:}]}%
\usepackage{array}

\newcommand{\x}{\mathbf{x}}

\newcommand{\reffig}[1]{Fig. \ref{#1}}
\usepackage{subfigure}

\PaperNumber{23-235}

\makeatletter
\renewcommand*\env@matrix[1][\arraystretch]{%
  \edef\arraystretch{#1}%
  \hskip -\arraycolsep
  \let\@ifnextchar\new@ifnextchar
  \array{*\c@MaxMatrixCols c}}
\makeatother

\begin{document}

\pagenumbering{arabic}
\title{Stochastic analysis of impulsive thrust uncertainties in the CR3BP}


\author{Sharad Sharan\thanks{Graduate student, Department of Aerospace Engineering, Pennsylvania State University, State College, PA-16802.}, 
Amit Jain\thanks{Graduate student, Department of Aerospace Engineering, Pennsylvania State University, State College, PA-16802.}, 
Roshan T. Eapen\thanks{Assitant Professor,  Department of Aerospace Engineering, Pennsylvania State University, State College, PA-16802.},
Puneet Singla\thanks{Professor, AIAA Associate Fellow, AAS Fellow, Department of Aerospace Engineering, Pennsylvania State University, State College, PA-16802.}, 
Robert G. Melton\thanks{Professor, AAS Fellow, AIAA Associate Fellow, Department of Aerospace Engineering, Pennsylvania State University, State College, PA-16802.}
}

\maketitle{} 		

\begin{abstract}
	This paper employs an alternate dynamical model of the circular restricted three body problem to quantify uncertainties associated with spacecraft thrusting maneuvers. A non-product quadrature scheme known as Conjugate Unscented Transform (CUT) is employed to determine the higher order system sensitivities through a computationally efficient data driven approach. Moreover, the CUT scheme, in conjunction with a sparse approximation method, is used to find an analytical representation of the time evolution of the state probability density function (pdf).



\end{abstract}

\section{Introduction}
	
    Several years of research and planning culminated into the commencement of the Artemis missions in the year 2022. With the launch of CAPSTONE and Artemis 1, the torch has been lit for human beings to return to the Moon, this time to stay. This indicates that cislunar space would soon be occupied by many spacecraft traversing it. A significant challenge in the volume of space collectively called cislunar, is that it is tedious to track spacecraft there. After every successful observation, it is required to propagate the uncertainties associated with the state of the spacecraft, so that tracking mechanisms are better prepared to make an observation at a later time. Should there be any off-nominal event experienced by the spacecraft between observations, then a tedious problem on a good day suddenly becomes intractable. Therefore, it is imperative to focus on techniques to quantify and propagate uncertainties associated with cislunar trajectories. Propagating the statistical moments of a distribution offers significant insight into the nature of a distribution. These insights are valuable for tracking strategies. This paper focuses on obtaining such insights by studying different kinds of uncertainties associated with an impulsive thrusting maneuver in cislunar space.    
    
    \par The Circular Restricted Three Body Problem (CR3BP) is a good starting model that is widely utilized in the study of cislunar trajectories. The CR3BP is a chaotic nonlinear system. Traditional schemes of propagating the central statistical moments rely heavily on linearized dynamics \cite{frueh2021cislunar}. These schemes are not ideal for a system like the CR3BP. In such cases, a widely used method of uncertainty propagation is the computationally intensive Monte Carlo (MC) method, where tens of thousands of simulations are run for each study \cite{vittaldev2016spacecraft, sabol2011probability}. There are several other techniques in literature that use alternate representations of uncertainty, such as the Gaussian mixture approach \cite{demars2013entropy,vishwajeet2014nonlinear, demars2014collision} and polynomial chaos expansions \cite{terejanu2008uncertainty, vittaldev2016spacecraft}. Methods like the State Transition Tensor Series (STTS) \cite{younes2012high} were proposed to replace the MC integrations to generate data for the propagation of moments. However, applying the STTS scheme involves the computation of higher order sensitivities, which is not a trivial process. In this paper, the Conjugate Unscented Transform (CUT) \cite{adurthi2015conjugate} technique is employed to propagate the central moments. Moreover, CUT is shown to compute the higher order sensitivities in a computationally efficient manner, using which the state of each sample in a distribution at any given time can be computed with ease.    
	\par Another factor in uncertainty propagation is the way in which the state of the spacecraft is defined \cite{frueh2021cislunar}. To this effect, a judicious choice of the coordinate system can alleviate many of the challenges in classical uncertainty propagation. A comprehensive theoretical survey of approaches to modeling and representing the restricted three body problem is provided by Szebehely \cite{szebehely2012theory}. However, the Cartesian coordinates are the widely used ones to model dynamics in the CR3BP. There has been very minimal numerical exploration using alternative coordinate systems. To this end, our previous work explored a modified curvilinear coordinate system consisting of three spherical position coordinates and two velocity pointing angles, to quantify uncertainty in the spacecraft's state from an alternate perspective \cite{sharan2022coordinate}. This curvilinear system, termed the Spherical-Velocity Angles Model (S-VAM), was inspired by the reduced order systems elaborated by Szebehely \cite{szebehely2012theory}.
	\par Propagating an initial Gaussian distribution of the states through the Cartesian CR3BP equations of motion generally leads to a spherical or ellipsoidal distribution. This implies that the uncertainty in velocity vector is defined in all directions. Having the velocity pointing angles themselves as state variables, the uncertainty in velocity vector direction can be defined within a cone, which makes more sense for applications where the general direction of heading of a spacecraft is known to lie within a cone of uncertainty, rather than having it be defined along all directions in a sphere of uncertainty. The S-VAM offers this capability. The advantage of such a perspective was outlined in our preceding work \cite{sharan2022coordinate}, with special focus on a data driven approach to compute the higher order state transition tensors of a trajectory. Using the results of that study, statistical moments up to the second order were computed. This paper follows up on that by investigating methodologies to propagate the state probability density function (pdf) through the CR3BP dynamics using the S-VAM, in addition to computing the higher order statistical moments using an optimal sampling of points provided by the CUT scheme. 

Sparse-based collocation methods have shown promise in capturing pdf of nonlinear systems \cite{mercurio2016conjugate,jain2023reachability,jain2023UP}. Such strategies are explored in this work with some modifications to determine the time-varying pdf of a distribution in the CR3BP. The interesting aspect of this lies in the fact that we end up with analytical approximations for the said pdf, with a minimal sampling of the collocation points. In the sparse-collocation method, the log-pdf is approximated using a polynomial-based expansion. Then a set of collocation points following the CUT scheme are generated to represent the pdf domain. A minimal number of basis functions is desired for our polynomial approximation. Approximating using a polynomial is nothing but fitting a surface to the multidimensional log-pdf. With a high degree of accuracy desired, and a high dimensional system, the number of combinations of basis functions becomes very large, thereby leading to overfitting issues. This work utilizes a sequential $l_{1}$-norm minimization routine to optimally select the dominant basis functions from the extensive dictionary of all combinations of basis functions.  

\par The organization of the paper is as follows. First, a brief overview of the CR3BP is provided, following which the alternate dynamical model, the S-VAM, from the preceding work \cite{sharan2022coordinate} is revisited. Subsequently, the numerical methodologies utilized in this paper are outlined. Following this, the results of applying the aforementioned methodologies to test scenarios in cislunar space is discussed. 

\section{Circular Restricted Three Body Problem}
\begin{figure}
\centering
\subfigure[The synodic frame]{\includegraphics[width=2.8in]{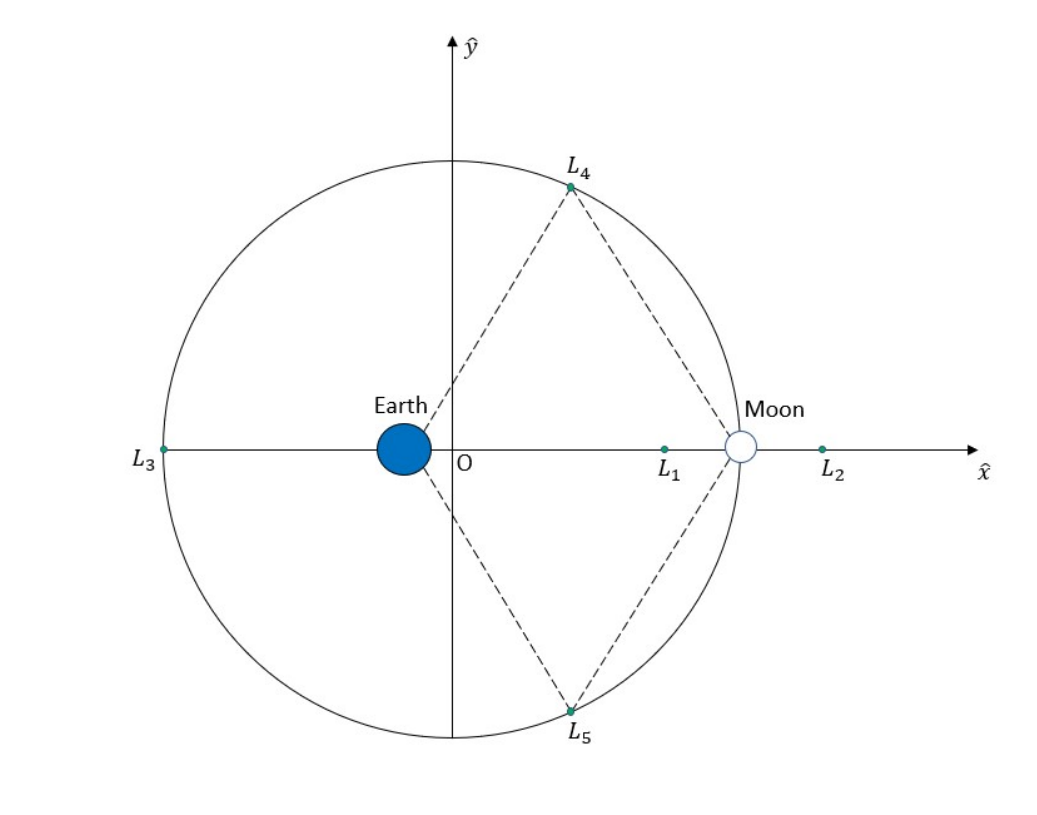}\label{fig:synodic}}
\subfigure[Velocity vector pointing angles of the S-VAM]{\includegraphics[width=2.5in]{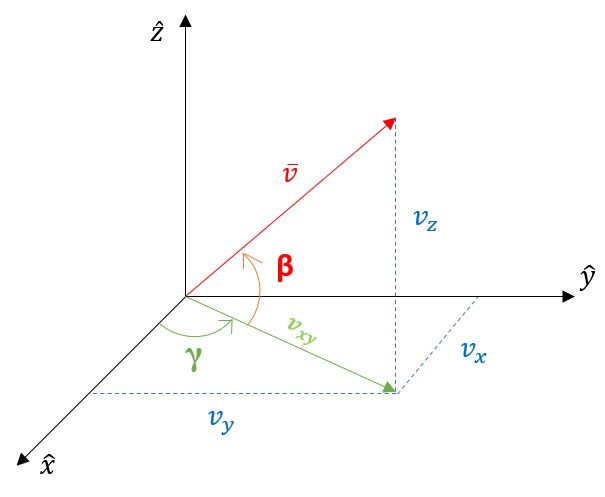}\label{fig:vpa_sketch}}
\caption{CR3BP schematic}\label{fig:CR3BP}
\end{figure}

The CR3BP is formulated in a frame that rotates along with the primaries, called the synodic frame. The Earth-Moon synodic frame is illustrated in Figure \ref{fig:synodic}. The $\hat{x}$ basis vector points from the origin, which is at the barycenter of the Earth-Moon system, toward the Moon. The $\hat{y}$ basis vector is perpendicular to it and lies in the plane of motion of the primaries as shown in Figure \ref{fig:synodic}. The $\hat{z}$ vector is given by the cross product of $\hat{x}$ and $\hat{y}$. 

A canonical system of units is employed where one length unit (LU) is equal to the distance between the two primaries and one time unit (TU) is chosen such that the mean motion of the primaries is unity. For the Earth-Moon system, $\SI{1}{LU} = \SI{384400}{\km}$ and $\SI{1}{TU} = \SI{4.3424}{}$ days. A pseudopotential is defined in the system as
%
\begin{align}
\Omega &= \frac{1}{2}(x^2 + y^2) + \frac{(1 - \mu)}{r_1} + \frac{\mu}{r_2}
\end{align}
where $r_1 = \sqrt{(x + \mu)^2 + y^2 + z^2}$, $r_2 = \sqrt{(x + \mu -1)^2 + y^2 + z^2}$, and $\mu = \frac{m_2}{m_1 + m_2}$. $r_1$ and $r_2$ are the magnitudes of the position vectors of the spacecraft relative to the Earth and the Moon, respectively. $m_1$ and $m_2$ are the masses of the Earth and the Moon, respectively. $\mu$ is the characteristic parameter in the synodic frame. For the Earth-Moon system, $\mu = 0.012151$. There exists an integral of motion in the CR3BP known as the Jacobi integral or the Jacobi constant ($C$).
\begin{align}
C = 2 \Omega - (v_x^2 + v_y^2 + v_z^2) \label{eq:C}
\end{align}
The Jacobi constant is the sole constant of motion in the CR3BP. Notice that all quantities given above are in terms of Cartesian coordinates. This is the traditionally used coordinate system to represent CR3BP dynamics.  However, in the process of numerical propagation, the Cartesian model accumulates errors over time, such that $C$ is no longer strictly maintained constant. The S-VAM prescribes $C$ explicitly in the equations of motion, thereby strictly holding it constant throughout the numerical propagation. The S-VAM equations are \cite{sharan2022coordinate}
\begin{align}
	\dot{r} &= \sqrt{2 \Omega - C} \left[ \cos{\phi} \cos{\beta} \cos{(\gamma - \theta)} + \sin{\phi} \sin{\beta} \right] \label{eq:rdot}\\[10pt]
	\dot{\theta} &= \frac{\sqrt{2 \Omega - C}}{r \cos\phi} \cos\beta \sin(\gamma - \theta) \label{eq:thetadot} \\[10pt]
	\dot{\phi} &= \frac{\sqrt{2 \Omega - C}}{r} \left[ \sin\beta \cos\phi - \sin\phi \cos\beta \cos(\gamma - \theta) \right] \label{eq:phidot}\\[10pt]
	\dot{\gamma} &=  \frac{1}{\sqrt{2 \Omega - C} \cos\beta} \left\lbrace r \cos\phi \sin\theta \cos\gamma \left( 1 - \frac{(1 - \mu)}{d_1^3} - \frac{\mu}{d_2^3} \right) \right.\nonumber\\ &- \left. \sin\gamma \left[ r \cos\phi \cos\theta - \frac{(1 - \mu)(r \cos\phi \cos\theta + \mu)}{d_1^3} - \frac{\mu(r \cos\phi \cos\theta - 1 + \mu)}{d_2^3} \right] \right\rbrace - 2 \label{eq:gammadotsph}\\[10pt]
	\dot{\beta} &= \frac{1}{\sqrt{2 \Omega - C}} \left\lbrace -r \sin\phi \cos\beta \left( \frac{(1 - \mu)}{d_1^3} + \frac{\mu}{d_2^3} \right) \right. \nonumber\\ &- \cos\gamma \sin\beta \left[ r \cos\phi \cos\theta - \frac{(1 - \mu)(r \cos\phi \cos\theta + \mu)}{d_1^3} - \frac{\mu(r \cos\phi \cos\theta - 1 + \mu)}{d_2^3} \right] \nonumber\\ &- \left. r \cos\phi \sin\theta \sin\gamma \sin\beta \left(  1 - \frac{(1 - \mu)}{d_1^3} - \frac{\mu}{d_2^3} \right) \right\rbrace \label{eq:betadotsph}
\end{align}
where, $d_1 = \sqrt{r^2 + \mu^2 + 2 \mu r \cos\phi \cos\theta}$ and $d_2 = \sqrt{r^2 + (\mu - 1)^2 + 2 r \cos\phi \cos\theta (\mu - 1)}$.
\par Eqs.\eqref{eq:rdot} - \eqref{eq:betadotsph} can now be numerically integrated with a guarantee that $C$ is preserved throughout the propagation. Using this model, a preliminary stochastic study associated with an impulsive thrust transfer from Low Earth Orbit (LEO) to an $L_2$ halo orbit was carried out previously \cite{sharan2022coordinate}. This work continues along this line and dives deeper into the relevant methodologies to compute metrics that allow one to obtain appropriate information useful for sensor tasking. The following section discusses the numerical methodologies employed in this paper to meet the aforementioned objective.
\section{Numerical methodology}
    The CUT technique is used extensively in this work for the following purposes: (i) Evaluation of statistical moments, (ii) computation of sensitivities, and (iii) propagation of state pdf.

CUT is an extended form of the unscented transform (UT), with modifications in the definition of the traditional sigma points. UT and CUT fall under a category of non-product quadrature schemes. A brief description of quadrature schemes, followed by the motivation for using CUT is given in the following section. 

\subsection{Evaluation of statistical moments}
The fundamental idea behind any quadrature scheme is that the integral of a function can be expressed as a weighted sum of function evaluations at specific points. These points and weights are what differentiate the various quadrature methods. Applying this idea to find the expectation value of $f(x)$,
\begin{align} \label{eq:expectation_quadrature}
	E[f(x)] &= \sum_{i=1}^{N} w_i f(x_i)
\end{align}        
where $x_i$ are the specific points where the function is evaluated, $w_i$ are the corresponding weights assigned and $N$ is the total number of points. Note that in a multidimensional system, $f$ and $x$ are both vectors. Expanding $f(x)$ into a Taylor series of $m^{th}$ order about a nominal solution ($x^*$), and equating the coefficients of the partials on both sides of the equation, one obtains a set of Moment Constraint Equations (MCEs).
%
\begin{align} \label{eq:MCEs}
	E[\delta x_{\alpha_1}] &= \sum_{i=1}^{N} w_i \delta x^{(i)}_{\alpha_1} \\
	E[\delta x_{\alpha_1} \delta x_{\alpha_2}] &= \sum_{i=1}^{N} w_i \delta x^{(i)}_{\alpha_1} \delta x^{(i)}_{\alpha_2} \\
	&\vdots \nonumber \\
	E[\delta x_{\alpha_1} \dots \delta x_{\alpha_m}] &= \sum_{i=1}^{N} w_i \delta x^{(i)}_{\alpha_1} \dots \delta x^{(i)}_{\alpha_m} \label{eq:HOMC}
\end{align}

Given the knowledge of the expectation values of $x$, the moments of $f(x)$ can be evaluated up to $m^{th}$ order accuracy by solving for a specific set of points and corresponding weights, such that the MCEs are satisfied \cite{hall2018probabilisticMS}.

    \par The traditional method of finding statistical moments is the Monte Carlo (MC) method. Several points from an initial probability density function are selected and propagated through the nonlinear system dynamics to evaluate the statistical moments at the final time. This scheme of random sampling means that the MCEs are approximately satisfied as the number of sampled points increases, but never exactly satisfied. For this reason, there are several deterministic methods that are used effectively. A brief outline of the deterministic approaches employed in determining the aforementioned points and weights can be found in literature \cite{hall2018probabilisticMS, sharan2022coordinate, jain2020computationally}. The well renowned ones include Gaussian quadrature \cite{golub1969calculation}, sparse grid methods \cite{gerstner1998numerical, heiss2008likelihood, jia2012sparse} and the UT \cite{julier2000new}. 
    
    \par Out of the methods listed above, only the UT allows us to define sigma points without the need to compute tensor products. This fact constraints the UT to be ineffective for higher dimensions. The sigma points lie symmetrically on the principal axes of the multidimensional variables, due to which all odd order moments are automatically satisfied. However, moving to cross dimensional moments greater than the third order, these sigma points are no longer sufficient to satisfy the MCEs. Thus, the CUT scheme is used in this work to evaluate multidimensional higher order moments. 
    
	\par The CUT is similar to UT in that it has points placed symmetrically on the principal axes such that odd-order moments are automatically satisfied. In addition to these axes, a conjugate set of axes is also defined, which contain several symmetrically placed points as well. These points along the conjugate axes help satisfy the cross-moment constraints in the case of multidimensional variables \cite{adurthi2012conjugate}. Altogether, the CUT points ensure the exact evaluation of multidimensional expectation integrals with a significantly smaller number of points than the previously discussed schemes. Adurthi et al. \cite{adurthi2012conjugate, adurthi2012conjugate1, adurthi2015conjugate1, adurthi2018conjugate} provide an extensive study with applications outlining the advantages of the CUT methodology over conventional quadrature rules. 

    \par The CUT scheme discussed here provides an optimal sampling of points to facilitate the computation of a trajectory's higher order sensitivities through a data driven approach. This approach is detailed in the following section.
	
\subsection{Computation of sensitivities}
A function $f(x)$ can be approximated by means of a Taylor series up to order $m$ about a reference solution$f(x^*)$ as 
\begin{align} \label{eq:Taylor}
	f(x) &= f(x^*) + \frac{\partial f(x^*)}{\partial x_{\alpha_1}}\delta x_{\alpha_1} + \ldots \frac{1}{m!} \frac{\partial^m f(x^*)}{\partial x_{\alpha_1} \ldots \partial x_{\alpha_m}}\delta x_{\alpha_1} \dots \delta x_{\alpha_m}
\end{align}
This can be succintly written as
\begin{align} \label{eq:approx}
	\vb{f}(\vb{x}) &\approx D \vb{\Phi}(\vb{x})
\end{align}  
where $D$ is a matrix of coefficients corresponding to the partial derivatives, and $\Phi(x)$ is a vector of polynomial basis functions. In general, evaluation of the partial derivatives to build the $D$ matrix is a computationally expensive process. However, if $f(x)$ and $\Phi(x)$ were known, a simple least squares (LS) procedure can be used to evaluate $D$. LS is a statistical approach and one requires data to initiate the process. It was earlier established that the CUT points were adept at satisfying higher order moment constraints. This essentially means that the CUT points, although few in number, are placed systematically in a domain of interest, such that the characteristics of the entire domain are well captured. Hence, the CUT points are chosen as the optimal sampling points at which to evaluate $f(x)$ and $\Phi(x)$, in order to generate the requied data for the LS procedure. From Eq. \eqref{eq:approx} the approximation error is given by
\begin{align}
	e_j &= f_j(x) - d_{ji} \Phi_i(x) \\
	\intertext{The cost function to minimize is}
	J &= \frac{1}{2}<e_j, e_j>
\end{align}
The inner product is with respect to a pdf that represents the domain of interest. Minimizing $J$,
\begin{align}
	\frac{\partial J}{\partial d_{jk}} &= 0 \\
	\Rightarrow \; 
	 <f_j(x), \, \Phi_k(x)> &= d_{ji} <\Phi_i(x), \, \Phi_k(x)> \\
	 d_{ji} &= \frac{<f_j(x), \, \Phi_k(x)>}{<\Phi_i(x), \, \Phi_k(x)>} \label{eq:Dmatrix}
\intertext{However, given a set of orthogonal polynomial basis functions, we can say that}
	<\Phi_i(x), \, \Phi_k(x)> &= 0, \; \text{for} \; i \neq k
\end{align}
This further reduces computational complexity, as the matrix of the inner products of these orthogonal polynomials (the normal matrix) is diagonal. Computing the inverse of a diagonal matrix is trivial, therefore, basis functions are preferred to be orthogonal. In addition to this, the normal matrix can be computed offline, thereby reducing on-the-fly computational load. In this work, the orthogonal polynomial basis functions up to degree four are considered. As a result of making $\phi_i(x)$ orthogonal, the coefficients computed are not exactly a reflection of the partial derivatives that appear in the Taylor series, nevertheless, they carry the same notion of sensitivities in the new basis system. Thus, using Eq. \eqref{eq:Dmatrix}, the sensitivity matrix $D$ is constructed.

\par The matrix $D$ can now be used to approximate the solution at any point in the vicinity of the reference trajectory as $D \Phi(\vb{\zeta})$, where $\vb{\zeta}$ is the state vector normalized in accordance with the distribution of CUT points. This normalization is done to alleviate numerical problems. In this way, the solution at the final time for all the samples in a distribution can be found using the sensitivity coefficients, henceforth referred to as the CUT-STTs. They are extremely useful, as they allow a significantly quicker function evaluation, as opposed to a series of integrations using an MC scheme.

\par Computation of the sensitivities is sufficient to commence the process of uncertainty propagation in a test scenario. These sensitivities allow one to approximate the states of all samples in a distribution using a polynomial series, thereby foregoing the need to integrate each sample individually. Using the integrated states, state histograms at the final time can be constructed to obtain insight into the consequence of the assumed initial uncertainties. However, using the CUT scheme, we can go one step further and develop methods to come up with analytical approximations of multidimensional pdfs. Such methods are discussed in the following section. 

\subsection{Propagation of state pdf}

State dynamics of an arbitrary system can be written as
\begin{equation}\label{DynamicalSystem_Cont}
    \dot \x = \mathbf{f}(\x,t)
\end{equation} 
Using the system flow $\vb{F}$, one can write
\begin{equation}\label{DynamicalSystem_Discrete}
    \mathbf{x}_{k} = \mathbf{F}( \x_0, t_{k}) 
\end{equation} 
Given the initial pdf $p(\x_0)$, the pdf $p(\mathbf{x}_k)$ at time $t_{k}$ can be evaluated as \cite{jazwinski2007stochastic}
\begin{equation}
p(\mathbf{x}_k) =p\left[\mathbf{x}_0=\mathbf{F}^{-1}(\mathbf{x}_k)\right]\left|\frac{\partial \mathbf{F}^{-1}}{\partial \mathbf{x}_k}\right|\label{eq:truePdf}
\end{equation}
where $\vb{F}$ is assumed to be an invertible, continuously differentiable mapping. 

The objective now is to develop a numerical approximation to Eq. \eqref{eq:truePdf}. In order to ensure smoothness, the structure of the state pdf is assumed to be
\begin{equation}
p\left(\mathbf{x}_{k+1}, t_{k+1}\right)=\delta p\left(\mathbf{x}_{k+1}, t_{k+1}\right) p\left(\mathbf{x}_{k+1}, t_{k}\right)
\end{equation}
where $p\left(\mathbf{x}_{k+1}, t_{k+1}\right)$ is the pdf of state $\x_{k+1}$ at time $t_{k+1}$, $p\left(\mathbf{x}_{k+1}, t_{k}\right)$ is the pdf of state $\x_{k+1}$ at time $t_{k}$, and $\delta p\left(\mathbf{x}_{k+1}, t_{k+1}\right)$ is the departure pdf of state $\x_{k+1}$ from time $t_{k}$ to $t_{k+1}$. An important condition to consider while approximating the pdf is that it must always be positive. Therefore, an exponential form of the pdf is assumed.
\begin{align}
    p(\vb{x}_{k+1}, t_{k+1}) &= e^{B(\vb{x}_{k+1}, t_{k+1})} \label{eq:pdf_exp}
\end{align}
Now, a polynomial approximation for the exponent $B(\vb{x}_{k+1}, t_{k+1})$ is sought. 
\begin{align}
B\left(\mathbf{x}_{k+1}, t_{k+1}\right) &= \delta B\left(\mathbf{x}_{k+1}, t_{k+1}\right)  + B\left(\mathbf{x}_{k+1}, t_{k}\right) \\
 B\left(\mathbf{x}_{k+1}, t_{k+1}\right) \approx \mathbf{c}_{k+1}^{T} \Phi\left(\mathbf{x}_{k+1} \right) &= \delta \mathbf{c}_{k+1}^{T} \Phi\left(\mathbf{x}_{k+1}\right)+\mathbf{c}_{k}^{T} \Phi\left(\mathbf{x}_{k+1}\right) \label{eq:del_ck}
\end{align}
where $\mathbf{c}_{k+1}$ is a vector of unknown coefficients at time $t_{k+1}$, $\delta\mathbf{c}_{k+1}$ is a vector of unknown departure coefficients from time $t_{k}$ to time $t_{k+1}$, $\mathbf{c}_{k}$ is the coefficients at time $t_{k}$, and $\Phi(\x_{k+1})$ is a vector of the chosen basis functions evaluated at state $\x_{k+1}$. The objective is to compute the departure coefficients $\delta\mathbf{c}_{k+1}$, which then serve as a weight function for the state pdf at the previous time. Thus, the smoothness constraint is enforced. Eq. \eqref{eq:del_ck} can be simplified to 
\begin{equation}\label{Eq:Ack+1=b}
\mathbf{A} \delta \mathbf{c}_{k+1}=\mathbf{b} 
\end{equation}
where,
\begin{equation}
\begin{gathered}
\mathbf{A}_{j}=\Phi^{T}\left(\mathbf{x}^j_{k+1}\right), \quad j=1,2, \ldots, N \\
\mathbf{b}_{j} = B(\vb{x}_{k+1}^j, t_{k+1}) -\Phi^{T}\left(\mathbf{x}^j_{k+1}\right)  \mathbf{c}_{k}, \quad j=1,2, \ldots, N
\end{gathered}
\end{equation}
The solution of Eq. \eqref{Eq:Ack+1=b} can be obtained using LS, which aims to find the best-fit solution for all the collocation points, referred to as the training data. The optimal value of coefficients $\delta \mathbf{c}_{{k+1}_{ls}}$ can be found by solving the weighted two-norm minimization problem
\begin{equation}
\begin{aligned}
\delta \mathbf{c}_{{k+1}_{ls}} = & \min_{\delta \mathbf{c}_{{k+1}}}  \left\|\mathbf{W}( \mathbf{A} \delta \mathbf{c}_{k+1} - \mathbf{b} ) \right\|_{2} \\
\end{aligned}
\end{equation}
where $\mathbf{W}$ is the weight matrix corresponding to the CUT points.



\par The LS solution solves for coefficients corresponding to all the basis functions in the dictionary. This leads to an overfitting issue, thereby introducing numerical inconsistencies. Moreover, since we are approximating the log pdf ($B$), the numerical errors become pronounced while finding the pdf as $e^B$. Therefore, a weighted $l_1$-norm minimization problem is posed to incorporate the minimum possible basis functions. Due to the equality constraint (Eq. \eqref{Eq:Ack+1=b}), this minimization problem is posed as
\begin{equation}
\begin{aligned}
&\min _{\delta \mathbf{c}_{k+1}}\left\|\mathbf{K} \delta \mathbf{c}_{k+1}\right\|_{1}\;\;
\text { subject to: } \;\;
& \left\|\mathbf{W}( \mathbf{A} \delta \mathbf{c}_{k+1} - \mathbf{b} ) \right\|_{2}  \leq \epsilon 
\end{aligned}
\end{equation}
where $\mathbf{K}$ is the coefficient weight matrix, and $\epsilon$ represents a threshold on the two-norm error of the sparse solution. The initial value of $\mathbf{K}$ can be set to 1 or based on any \textit{a-priori} knowledge gained from the least-squares solution. For subsequent iterations, $\mathbf{K}$ can be modified to penalize coefficients smaller than a predetermined threshold $\Delta_{s}$.
$$\mathbf{K}= \frac{1} {(\delta \mathbf{c}_{{k+1}_{ls}} + \eta)}$$ 
where $\eta$ is an arbitrary small number in order to prevent singularities. The value of $\epsilon$ is chosen to provide a flexible solution to the optimization problem by trading off approximation error with sparsity. Algorithm \ref{SparseSolution_algo} describes the procedure to compute a minimal polynomial approximation for the log-pdf. The user-defined parameter $\Delta_{s}$ represents the desired difference between the two norms of consecutive coefficients. The algorithm eventually arrives at an optimal solution $\mathbf{c}_{k+1}^*$, when the difference between the two norms of consecutive coefficient values, $\delta$, is less than $\Delta_{s}$.

\begin{algorithm}
\caption{Iterative weighted $l_1$-norm optimization \label{SparseSolution_algo}}
  \begin{algorithmic}[1] 
    \INPUT $\mathbf{A}, \mathbf{b},\mathbf{W},\mathbf{c}_{{k+1}_{ls}}, \Delta_{s}, \alpha, \epsilon , \eta$
    \OUTPUT $\mathbf{c}_{k+1}^*$
    \STATE \textbf{Initialization}   $\mathbf{K}= \frac{1} {(\mathbf{c}_{{k+1}_{ls}} + \eta)}, \delta=1  $ 
    \STATE  compute $\mathbf{c}_{k+1}^- =  \displaystyle \min _{\mathbf{c}_{k+1}}\left\|\mathbf{K} \mathbf{c}_{k+1}\right\|_{1} $  
    
     $\quad \qquad  \text{subject to:} \quad \left\|\mathbf{W}( \mathbf{A} \mathbf{c}_{k+1} - \mathbf{b} ) \right\|_{2}  \leq \epsilon$ 
    \WHILE{$\delta \ge \Delta_{s}$}
      \STATE Update $\mathbf{K}= \frac{1}{(\mathbf{c}_{k+1}^- + \eta)}$, to find $\mathbf{c}_{k+1}^+  =  \displaystyle \min _{\mathbf{c}_{k+1}}\left\|\mathbf{K} \mathbf{c}_{k+1}\right\|_{1} $  
      
      $\quad \quad\quad \quad\quad \qquad \quad \qquad\quad \qquad  \text{subject to:} \quad \left\|\mathbf{W}( \mathbf{A} \mathbf{c}_{k+1} - \mathbf{b} ) \right\|_{2}  \leq \epsilon$
      \STATE  Compute $\delta = \| \mathbf{c}_{k+1}^+ - \mathbf{c}_{k+1}^-\|_2 $
      \STATE $\mathbf{c}_{k+1}^- = \mathbf{c}_{k+1}^+ $ 
    \ENDWHILE
    \STATE $\mathbf{c}_{k+1}^* = \mathbf{c}_{k+1}^-$
  \end{algorithmic}
\end{algorithm}

After obtaining the sparse coefficients $\vb{c}_{k+1}^*$, it is possible to separate the dominant coefficients from the non-dominant coefficients by choosing a user-defined coefficient threshold $\delta_{rs}$. The non-dominant coefficients are those with magnitudes lesser than $\delta_{rs}$, and can be neglected for computational purposes. Therefore, a new minimal representation of the basis functions $ \vb{A}_{rs}\in \mathbb{R}^{N\times m_r} $ corresponding to the $m_r$ non-zero coefficients can be constructed. These are termed reduced sparse (RS) coefficients, and can be calculated using least-squares as
\begin{equation}
\begin{aligned}
\mathbf{c}_{{k+1}_{rs}} =& \mathbf{A}_{rs}^{\dagger}  \mathbf{b}
\end{aligned}
\end{equation}
where $\mathbf{A}_{rs}^{\dagger}$ is the pseudo-inverse. Now, these coefficients along with the minimal basis functions yield the required polynomial approximation of the state pdf.
\par The $l_1$-norm minimization problem addressed in Algorithm \ref{SparseSolution_algo} is convex, and must be invoked at each time instant to find the time-varying state pdf. Various convex optimization solvers like Sedumi, SDPT3, Gurobi, MOSEK, and GLPK, can be used  to solve this problem\cite{boyd2004convex,cvx}. While finding a sparse solution, it must be noted that a single solver may not be able to provide the optimal solution, as each solver has different capabilities. Since numerical methods for convex optimization are not exact, the results are computed within a numerical precision or tolerance that has been predefined.
\par In the process detailed above, the CUT points are chosen as the collocation points where training data is obtained. As mentioned earlier, this is because the CUT points provide the minimal number of points required to appropriately capture the characteristics of the domain of interest.

\par The succeeding sections of the paper discuss the results of applying these methodologies to a couple of case studies in the CR3BP.


\section{Results and Discussion}
    Two different scenarios with different parameter uncertainties associated with an impulsive thrusting maneuver are studied in this paper. Note that scenario A was explored in the previous work \cite{sharan2022coordinate} predominantly for the accuracy of the CUT-STTs, however, the pdf propagation scheme and higher order moments computation was not investigated at all. Since these features are explored in this paper, a quick summary of scenario A is provided before moving on to the details of pdf propagation.  Three cases are considered in scenario A, which deals with uncertainty in the position of the spacecraft at impulse, and in the direction of the impulse. These cases have exaggerated values of uncertainty. Such values were deliberately assumed to explore the effectiveness of the CUT-STTs in approximating the propagation of each sample in the state.
    
    \par Subsequently, scenario B deals with plausible values of uncertainty in the direction of the impulse, its magnitude and the time of its occurrence. This scenario raises the bar for the CUT-STTs, as considerably fewer sampling points are chosen using the CUT scheme in comparison to scenario A, yet, good polynomial approximations of the states of all samples are generated, without the need for integrating any sample in the distribution. Following the testing of CUT-STTs in both scenarios, the propagation of state pdf is tested on Case 1 of scenario A. The particulars of these tests are discussed in detail in the following sections.
    
\subsection{Testing the CUT-STTs}
\begin{figure}[h]
\centering
\includegraphics[width=4.5in]{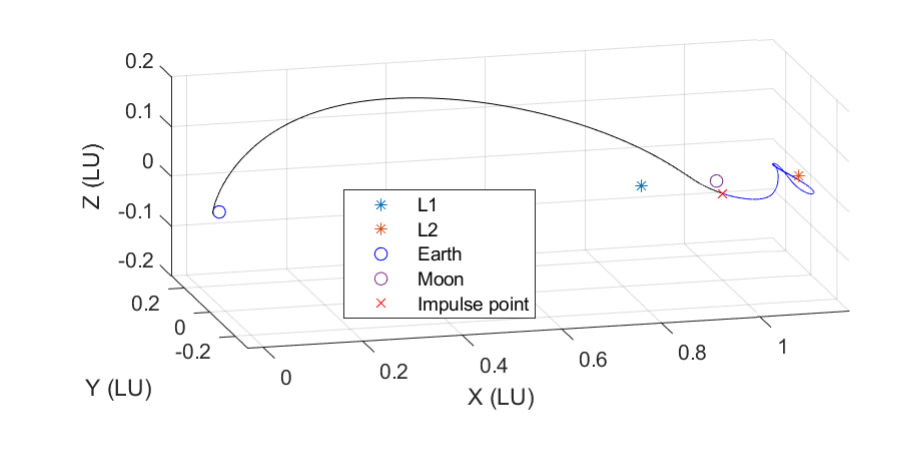}
\caption{Reference trajectory from LEO to an $L_2$ halo orbit}
\label{fig:refTraj}
\end{figure}

	A reference transfer from a low Earth orbit to an $L_2$ halo orbit is chosen, as illustrated in Fig. \ref{fig:refTraj}. An impulsive $\Delta V$ is applied to this translunar trajectory to effect insertion into a trajectory on the stable manifold of the halo orbit. Due to the nature of S-VAM, uncertainty in the magnitude of the velocity vector can be handled, independent of the uncertainty in its direction. The following study deals with the uncertainty associated with the firing angle of the thruster, and the position of the spacecraft at the impulse point. In order to capture this uncertainty, a uniform distribution of samples corresponding to $\pm [\SI{5}{\km}, \SI{0.2}{\degree}, \SI{0.2}{\degree}, \SI{5}{\degree}, \SI{5}{\degree}]$ in $[r, \theta, \phi, \gamma, \beta]$ about the impulse point on the reference trajectory is generated. These uncertainty bounds are deliberately chosen to be considerably higher than realistic navigation errors or attitude errors at the time of the impulse. This is to showcase the effectiveness of the aforementioned numerical methodologies applied from the S-VAM perspective. Note that this setup assumes that the energy imparted by the thruster is same as that of the reference transfer. 
    \par First, \num{1000000} Monte Carlo samples are integrated for a time period of two days using the S-VAM. This is followed by multiple sets of simulations, that are carried out by subsequently reducing the number of samples by an order of ten. 
\par
\begin{figure}[h]
        	\centering
        	\includegraphics[width=3.5in]{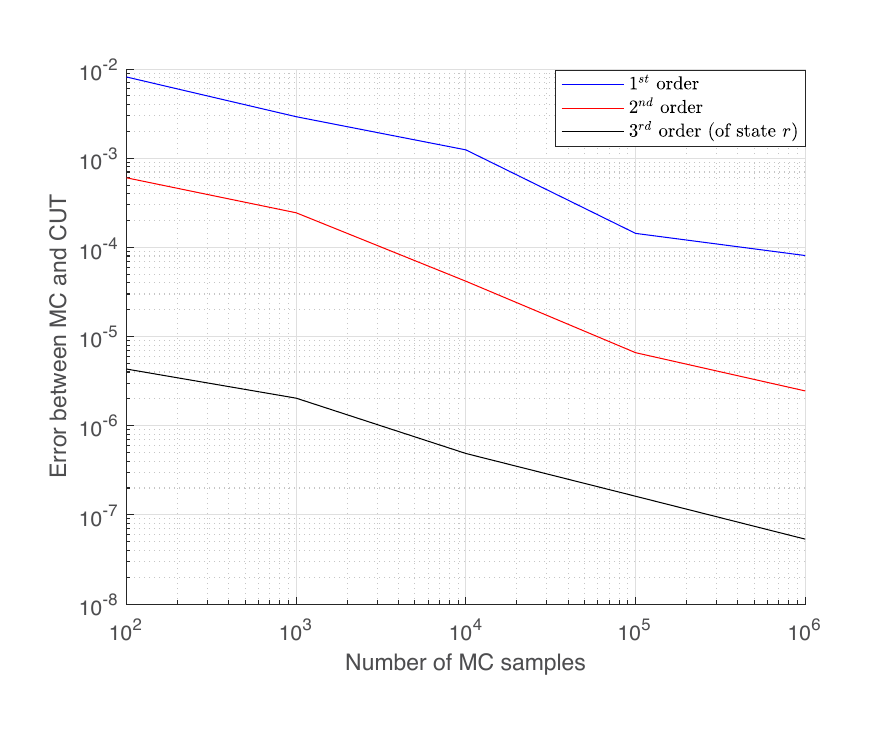}
        	\caption{Two-norm error between the CUT and MC computations of the centralized statistical moments up to third order}
        	\label{fig:statMomComp}
   \end{figure}

 Centralized statistical moments up to the fourth order are calculated for each set of MC simulations, and compared with those obtained using the CUT methodology. The two norm of the corresponding errors in these moments between the two methods are plotted as illustrated in Fig. \ref{fig:statMomComp}. For the third order moment, the matrix corresponding to the state $r$ is pulled from the error tensor, and its two norm is plotted. It can be observed that as the number of samples increases, the error between MC and CUT reduce. This is indicative of the fact that the CUT points enable an accurate calculation of the moments without the need for numerous samples like the MC method. On the other hand, with the MC scheme, greater the sample size, better is the estimate. Thus, the error reduces as the number of MC points increases. Visualizing the fourth order moment tensor becomes tricky, nevertheless, it was computed using the CUT and MC methods. These statistical moments provide information to tracking strategies for better situational awareness. In addition to the statistical moments, the final state of any random initial sample can be approximated using the CUT-STTs quickly and efficiently.   
\par
\begin{figure}[h]
        \centering
        \includegraphics[width=4in]{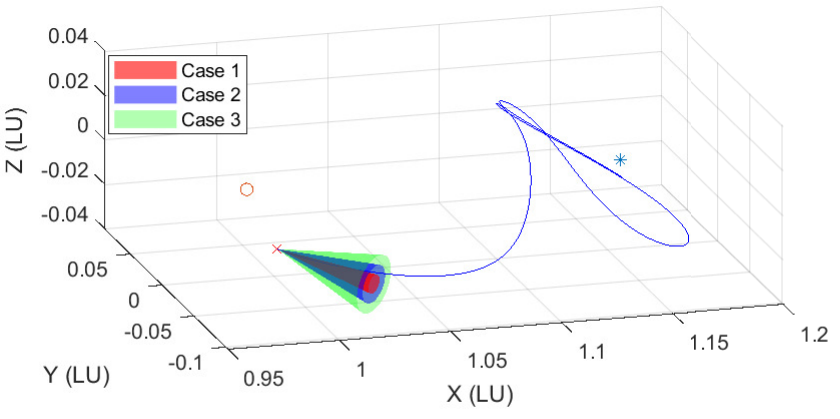}
        \caption{Cone of uncertainty ($\gamma$, $\beta$)}
        \label{fig:refTrajCones}
   \end{figure}

\begin{figure}[H]
   \centering
   \includegraphics[scale=0.75]{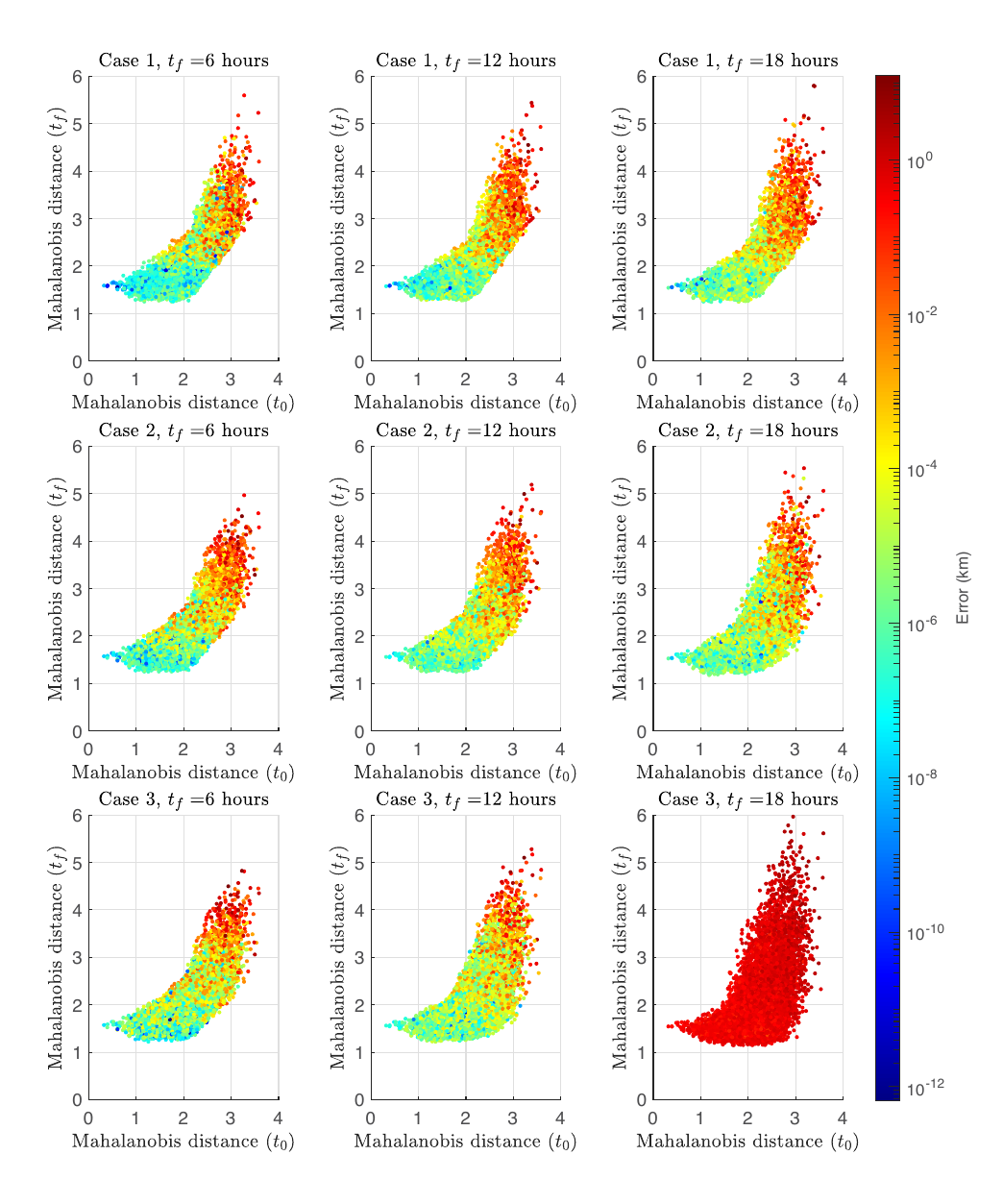}
   \caption{Position error (in km) between MC integration and CUT approximation plotted against the initial and final Mahalanobis distances of the samples}
   \label{fig:MDs}
\end{figure}

 In order to learn the accuracy of the CUT-STT approximation, the two norm of the error in position is plotted against the initial and final Mahalanobis distances of the samples, as shown in Fig. \ref{fig:MDs}. This is the error between the position computed by an MC integration and the position approximated using the CUT-STTs. The Mahalanobis distance ($M_d$) of each sample is calculated as 
\begin{align}
	M_d &= \sqrt{(\vec{x}_t - \vec{\mu}_t)^T \Sigma_t^{-1} (\vec{x}_t - \vec{\tilde{\mu}}_t)}
\end{align}
where $\mu_t$ is the mean vector of the state distribution at time $t$, and $\Sigma_t$ is the covariance of the distribution. $M_d$ is a measure of how far the multi-dimensional sample is from the mean. Three different cases are studied, with each case differing by the degree of uncertainty in the velocity pointing angles. The samples in cases $1$, $2$ and $3$ are associated with a maximum uncertainty of $\pm \SI{5}{\degree}, \pm \SI{10}{\degree}$ and $\pm \SI{15}{\degree}$, respectively, in $\gamma$ and $\beta$. This is illustrated in Fig. \ref{fig:refTrajCones} by means of the cones of uncertainty. It can be discerned from these cones whether or not a spacecraft could possibly reach the intended terminal manifold given a directional uncertainty in the impulsive maneuver. In this region of the CR3BP, such an analysis is important, because if the spacecraft has enough energy to reach the family of halos around $L_2$ and is not heading the right way, then it could possibly go through the $L_2$ gateway and out of the system. Uncertainty in the position of the spacecraft during the impulse is considered to be within $\pm[\SI{5}{\km}, \SI{0.2}{\degree}, \SI{0.2}{\degree}]$ about the reference impulse point for all the cases.
\par The time at the impulse point is designated $t_0$, and the time of propagation of the samples is designated $t_f$. $M_d$ is calculated for the samples at $t_0$ and $t_f$ and plotted along the $x$ and $y$ axes, respectively, in Fig. \ref{fig:MDs}. The colorbar is representative of the position error in km between the MC integration and the CUT approximation. From Fig. \ref{fig:MDs}, it can be observed that for cases $1$ and $2$, the CUT-STTs capture dynamics accurate to an order of \SI{e-1}{\km} up to $t_f = \SI{18}{\hour}$, even for samples that are at $5 \sigma$ at that $t_f$. Considering only up to $3 \sigma$ points for $t_f = \SI{12}{\hour}$, the maximum error across all cases is on the order of \SI{e-2}{\km}. These results are promising in terms of accuracy of the CUT-STTs.  

\par As the cone of uncertainty grows wider and $t_f$ increases, the accuracy of the CUT-STTs decreases, as illustrated in Fig. \ref{fig:MDs} for case $3$. One can conclude this from the sharp increase in error magnitude between $t_f = \SI{12}{\hour}$ and $t_f = \SI{18}{\hour}$ for case $3$. However, in a practical sense, an observation is likely to be made to update the position of the spacecraft within \SI{12}{\hour}.

\par
\begin{figure}[h]
   \centering
   \includegraphics[scale=0.85]{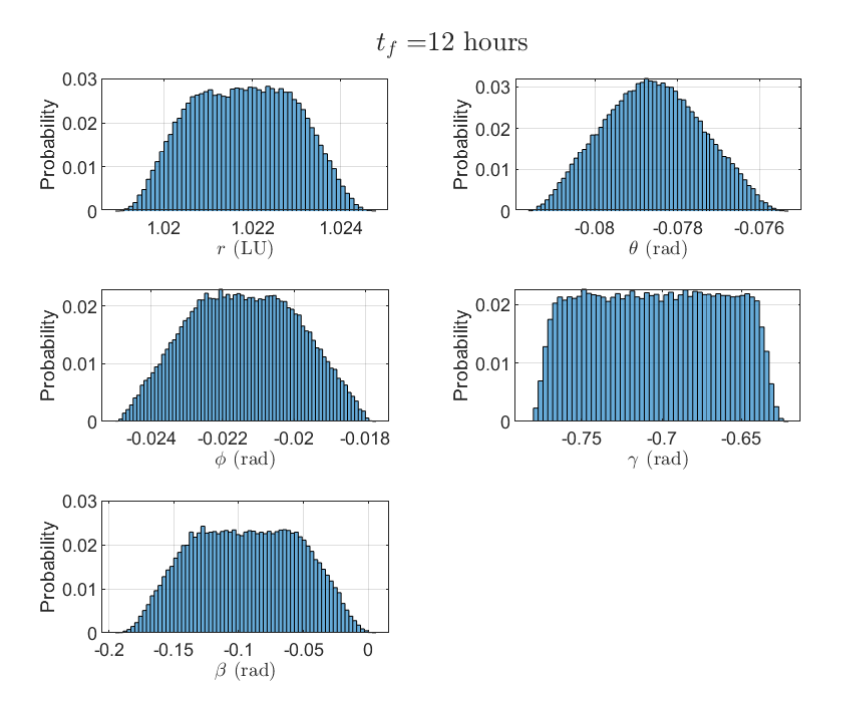}
   \caption{Histograms of the states at $t_f = \SI{12}{\hour}$ - Scenario A (case 1)}
   \label{fig:hist_A}
\end{figure}

 The analysis carried above sheds light on the duration and degree of accuracy of the CUT-STTs for the aforementioned cases. They can now be used to obtain stochastic insight that will aid tracking strategies. It would be beneficial to know the consequences of an off-nominal burn at some time in the future, so that the sensors can be tasked to re-acquire the spacecraft at the next instance of observation. In order to do that, useful information can be gathered from histograms of the states, which are illustrated in Fig. \ref{fig:hist_A} for case 1. 
\par It is to be noted that these histograms were plotted using states computed by the CUT-STTs. A total of \num{e5} samples were propagated using the CUT-STTs in a matter of seconds to obtain the histograms in Fig. \ref{fig:hist_A}. Greater the number of samples, better is the estimate of the distribution. However, with increased sample size comes great computational burden, if one relies on the traditional MC scheme.
 The CUT-STTs offer an efficient way to obtain insight by propagating a huge sample of states, without the computationally expensive aspect of it.  This can be realized by clocking the computation time for each method, as shown in Table \ref{tab:comTime}. All computations were run on an Intel Xeon E5-2695 processor and parallelized over 32 cores. There is no parallelization necessary while using the CUT scheme. Moreover, for this case study of a misaligned thruster, the S-VAM, being a reduced order model, offers further reduction in computation required by the CUT-STTs.
\begin{table}[h]
	\centering
        \caption{Comparison of computation time}
        \medskip
	\begin{tabular}{|c|c|c|}
		\hline
            \multirow{2}{*}{\textbf{Number of samples}}  
            & \multicolumn{2}{c|}{\textbf{Computation time} (sec)}  \\  \cline{2-3}
		  & \textbf{Monte Carlo}  & \textbf{CUT} \\
		\hline
		\num{1000000} & \num{507.06} & \num{9.43} \\
		\hline
		\num{100000} & \num{80.41} & \num{0.84} \\
		\hline
	\end{tabular}
	\label{tab:comTime}
\end{table}
\begin{table}[h]
	\centering
        \caption{Statistical moments evaluated using CUT8 points - Scenario A (case 1)}
	\begin{tabular}{|c|c|c|c|c|}
		\hline
		 & \textbf{Mean} & \textbf{Variance} & \textbf{Skewness} & \textbf{Kurtosis}\\
		\hline
		$\boldsymbol{r}$ & $1.021809$ & $\num{1.313916e-6}$ & $0.007993$ & $2.111832$\\
		\hline
		$\boldsymbol{\theta}$ & $-0.078584$ & $\num{1.373864e-6}$ & $0.046122$ & $2.382490$\\
		\hline
		$\boldsymbol{\phi}$ & $-0.021397$ & $\num{2.264103e-6}$ & $0.017949$ & $2.195804$\\
		\hline
            $\boldsymbol{\gamma}$ & $-0.703372$ & $\num{1.633346e-3}$ & $-0.004865$ & $1.825565$\\
            \hline
            $\boldsymbol{\beta}$ & $-0.094022$ & $\num{1.641397e-3}$ & $0.000463$ & $2.069970$\\
            \hline
	\end{tabular}
	\label{tab:statMomScenarioA}
\end{table}
\par The statistical moments corresponding to the final distribution of the states are given in Table \ref{tab:statMomScenarioA}. It is to be noted that these moments are computed using a set of \num{455} CUT points only, yet they capture the state distributions as shown in Fig. \ref{fig:hist_A} very well. In Table \ref{tab:statMomScenarioA}, the mean for each state is a measure that stands out as being readily comparable visually with the histogram of the corresponding state. Another measure that can be visually appreciated right away is the kurtosis of $\theta$. It can be seen from Fig. \ref{fig:hist_A} that the histogram of $\theta$ closely resembles a normal distribution. For normal distributions, kurtosis $= 3$. It can be seen from Table \ref{tab:statMomScenarioA} that the value of kurtosis for $\theta$ is the closest to three among all the states.

\par The covariance, skewness and kurtosis corresponding to the multivariate distribution were computed following the approach given in Eq. \eqref{eq:HOMC}. However, it is harder to visualize the full skewness and kurtosis tensors, thus, the univariate moments for each state are presented in Table \ref{tab:statMomScenarioA}.  

    Scenario A presented a very specific case of a misaligned impulse and studied its consequences for a particular reference trajectory. Primarily, it established the effectiveness of the CUT methodology and the usefulness of the S-VAM. The next scenario (scenario B) intends to study a wider array of uncertainties associated with a maneuver. For better contrast, the same reference trajectory from scenario A is used here. 
\begin{figure}[h]
   \centering
   \includegraphics[scale=0.75]{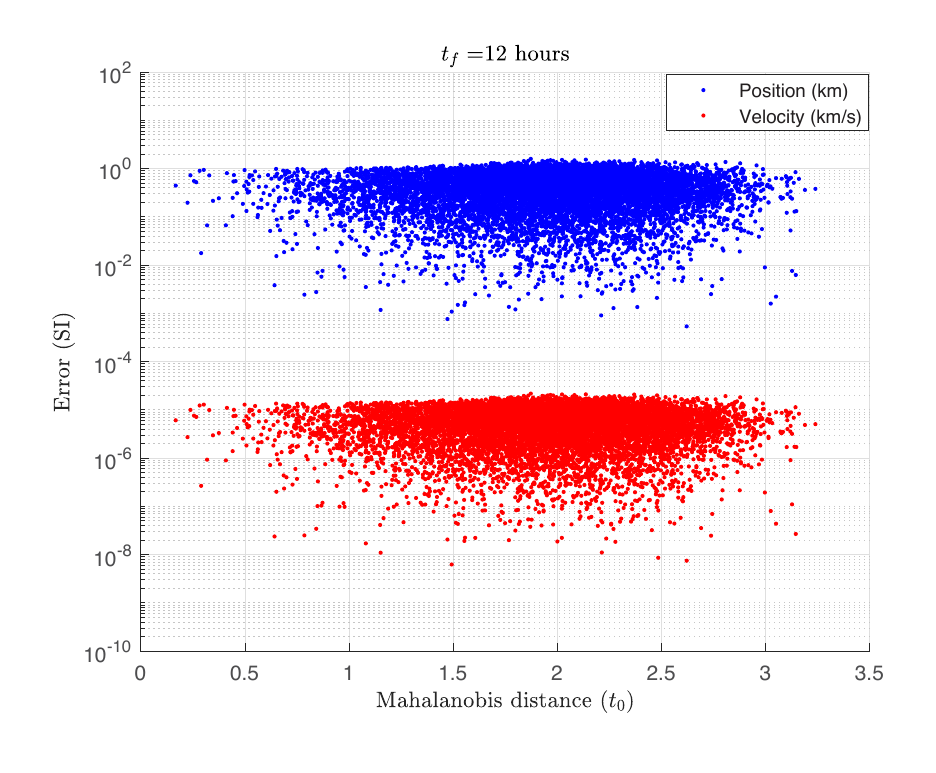}
   \caption{Position error (in km) between MC integration and CUT approximation - Scenario B}
   \label{fig:MDs_B}
\end{figure}
    \par The total thrust delivered by the spacecraft's engine is considered to be \SI{500}{\newton}, and the mass of the spacecraft at the start of the maneuver is assumed to be \SI{600}{\kg}. The reference trajectory includes a burn that delivers a nominal $\Delta V$ of \SI{677.5}{\meter \per \second} over a burn time of \SI{813}{\second}. Considering the magnitude of flight time involved in this cislunar trajectory, which is on the order of days, this maneuver can be approximated as being impulsive. Uncertainty is considered in the parameters $[t_1, t_b, \gamma, \beta]$, where $t_1$ is the time at which the maneuver occurs, and $t_b$ is the engine burn time. The uncertainty in burn time is to account for the incorrect magnitude of $\Delta V$ imparted. Uncertainties in $\gamma$ and $\beta$ account for the incorrect direction of the $\Delta V$. In this simulation, a maximum uncertainty in $t_1$, $t_b$, $\gamma$ and $\beta$ of $\pm \SI{60}{\sec}$, $\pm \SI{5}{\sec}$, $\pm \SI{2}{\degree}$, and $\pm \SI{2}{\degree}$ are considered about the nominal, respectively. $\num{100000}$ samples are generated within these uncertainty bounds and propagated for twelve hours. 

    \par The CUT-STTs computed in this scenario are sensitivities of the final states with respect to variations in the initial stochastic parameters, i.e., $[t_1, t_b, \gamma, \beta]$. In scenario A, the stochastic parameters were the initial states, therefore, the CUT-STTs obtained there were a reflection of the traditional state transition tensors. The primary advantage of these sensitivities with respect to the initial stochastic parameters, lies in the fact that the effect of any variation in the said parameters can be mapped to the final states rapidly.  
\begin{figure}[h]
   \centering
   \includegraphics[scale=0.8]{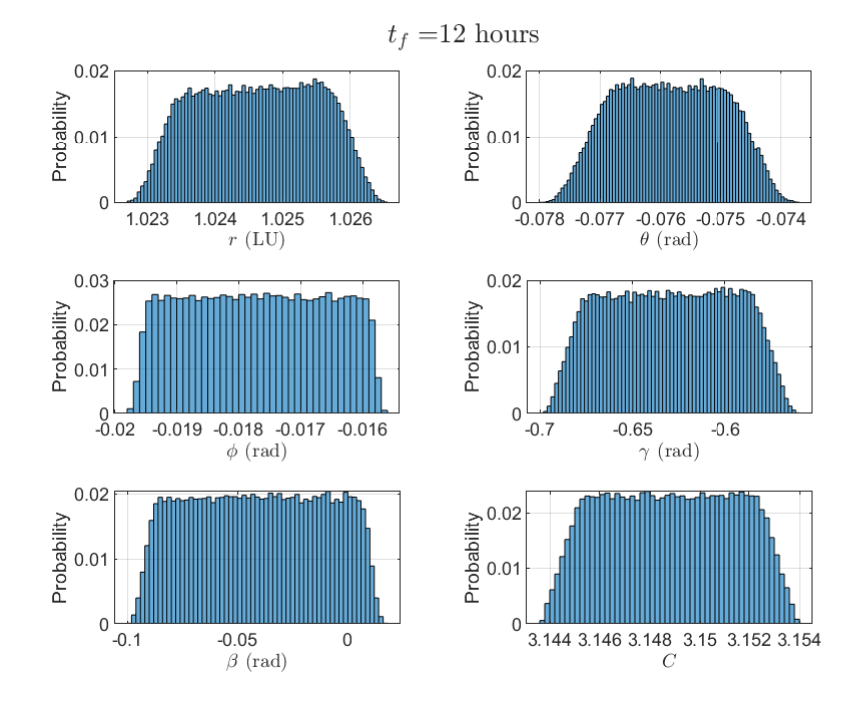}
   \caption{Histograms of the states at $t_f = \SI{12}{\hour}$ - Scenario B}
   \label{fig:hist_B}
\end{figure}
\begin{table}[h]
	\centering
        \caption{Statistical moments evaluated using CUT8 points - Scenario B}
	\begin{tabular}{|c|c|c|c|c|}
		\hline
		 & \textbf{Mean} & \textbf{Variance} & \textbf{Skewness} & \textbf{Kurtosis}\\
		\hline
		$\boldsymbol{r}$ & $1.024617$ & $\num{7.387324e-7}$ & $-0.028405$ & $1.922358$\\
		\hline
		$\boldsymbol{\theta}$ & $-0.075851$ & $\num{7.504202e-5}$ & $0.006545$ & $2.026064$\\
		\hline
		$\boldsymbol{\phi}$ & $-0.017698$ & $\num{1.213942e-6}$ & $-0.001075$ & $1.814514$\\
		\hline
            $\boldsymbol{\gamma}$ & $-0.630171$ & $\num{1.071206e-3}$ & $-0.011119$ & $1.874504$\\
            \hline
            $\boldsymbol{\beta}$ & $-0.040751$ & $\num{8.891200e-4}$ & $-0.007606$ & $1.826738$\\
            \hline
            $\boldsymbol{C}$ & $3.148817$ & $\num{6.570981e-6}$ & $-0.003734$ & $1.871113$\\
            \hline
	\end{tabular}
	\label{tab:statMomScenarioB}
\end{table}
    \par The samples are also propagated by MC integration in order to compare with the CUT-STT approximation. The errors in position and velocity of the samples at the final time are plotted against the corresponding initial Mahalanobis distances of the samples in Fig. \ref{fig:MDs_B}. The maximum error in velocity is on the order of \SI{e-5}{\km \per \s}, while the maximum position error is found to be two kilometers. 
     \par Considering the fact that the final distribution of the samples spans a distance that is several orders of magnitude greater than two kilometers, this error in position is relatively small, and does not affect the nature of the distribution in a significant way. Moreover, this is only the maximum error in position observed among all the samples. 
    
   \par The root mean squared error (RMSE) in position is found to be \SI{523}{\m}. This accuracy of the CUT-STT after $12$ hours of propagation is very promising, as the RMSE in position is only on the order of navigation errors that are usually experienced by spacecraft in cislunar space. Therefore, the effect of any initial uncertainty in the parameters ($[t_1, t_b, \gamma, \beta]$) can be rapidly mapped to states at a final time using the CUT-STTs by means of simple function evaluations, completely forgoing numerical integration. Using these states, histograms can be plotted as illustrated in Fig. \ref{fig:hist_B}. This proves better than the traditional way to obtain these histograms, which would be to propagate each sample through a numerical integrator, which in turn requires a considerably higher computational power. Although the sample size is $\num{100000}$ in the domain of interest, it takes only $161$ CUT points to capture the statistical moments, given in Table \ref{tab:statMomScenarioB}, thereby highlighting the computational efficiency of CUT. 
    
    \par The two scenarios discussed in this section showcased the effectiveness of the S-VAM in handling directional uncertainties, as well as the capabilities of the CUT-STTs in approximating the system dynamics. It was also observed that the optimal deterministic sampling provided by the CUT points was sufficient to evaluate statistical moments up to the fourth order with high accuracy that cannot be matched exactly by random sampling techniques.

\begin{figure}[h]
\subfigure[Magnitude of LS coefficients]{\includegraphics[width=0.5\textwidth]{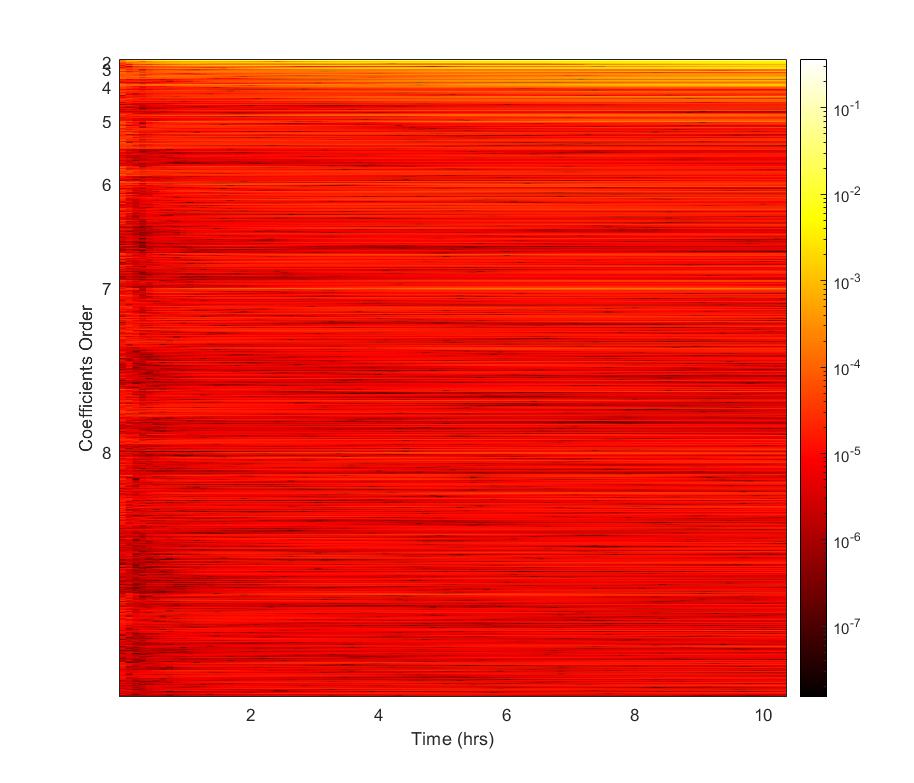}\label{fig:heatmap_LS}}
\subfigure[Magnitude of RS coefficients]{\includegraphics[width=0.5\textwidth]{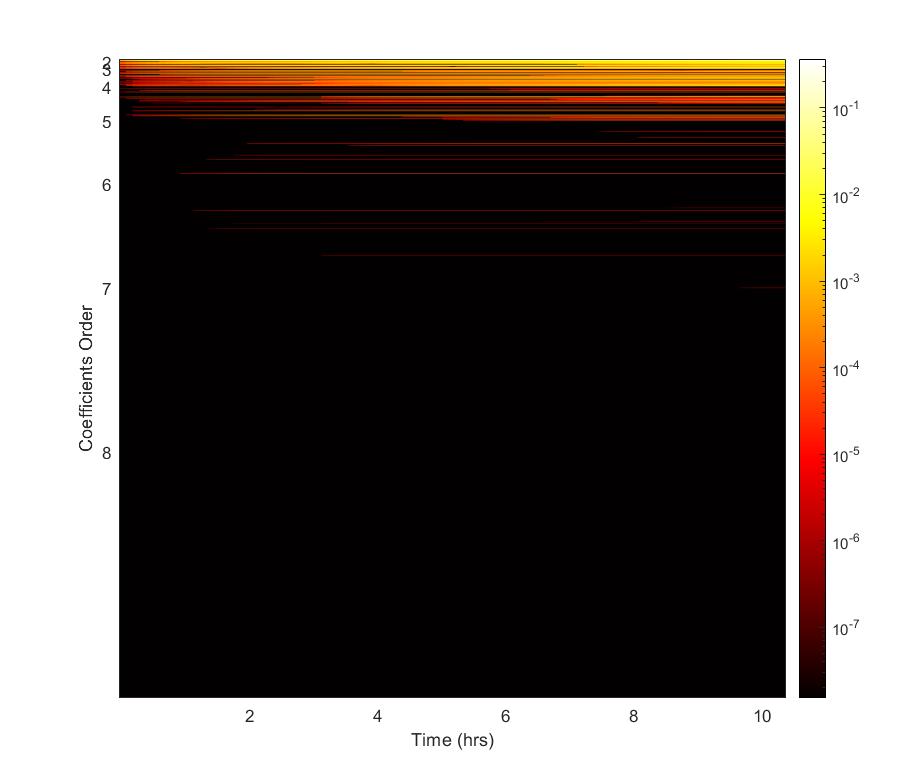}\label{fig:heatmap_RS}}
    \subfigure[Evolution of Non-Zero Coefficients]{\includegraphics[width=0.5\textwidth]{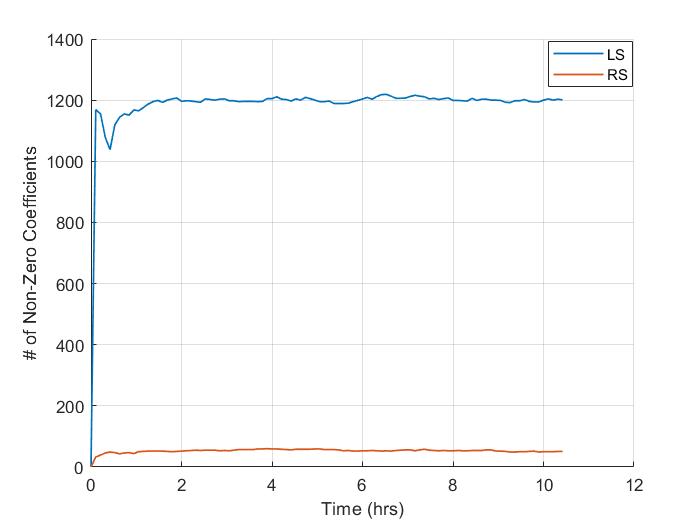}\label{fig:NNZ}}
    \subfigure[Comparing training and testing relative error]{\includegraphics[width=0.5\textwidth]{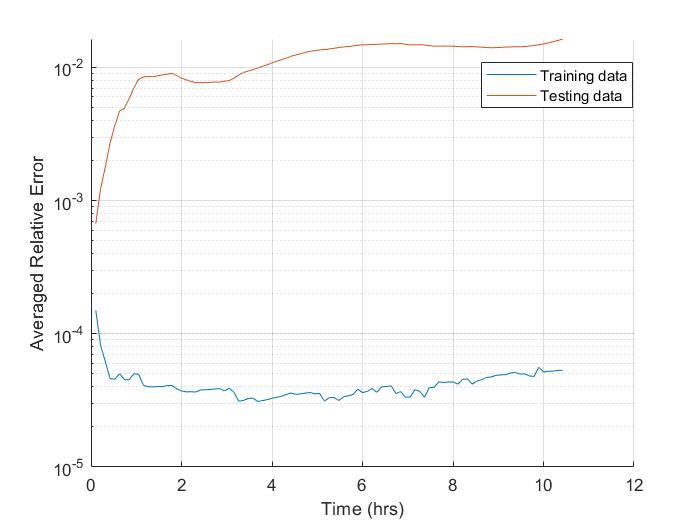}\label{fig:TestTrainError}}
\caption{Propagation of state pdf}\label{fig:NNZ_heatmap}
\end{figure}

\subsection{Testing state pdf propagation}
\par In order to test the methodology outlined for propagating the state pdf, case 1 of scenario A is chosen. Basis functions up to $8^{th}$ order are considered for the analytical approximation. The sparse algorithm parameters (as outlined in Algorithm \ref{SparseSolution_algo}) are taken to be $\epsilon = \num{1e-6}$, $\eta = \num{1e-5}$, $\Delta_{s}= \num{1e-5}$ and $\delta_{rs} = \num{1e-5}$.

\par The evolution of the LS and RS coefficients are illustrated using heatmaps in \reffig{fig:heatmap_LS} and \reffig{fig:heatmap_RS}, respectively. The colorbar is indicative of the magnitude of the coefficients. Essentially, it indicates how much an $n^{th}$ order basis function is participating in the analytical approximation. It can be seen from \reffig{fig:heatmap_LS} and \reffig{fig:heatmap_RS} that the LS solution utilizes most of the coefficients, while the RS solution only captures the dominant coefficients from the extensive dictionary of corresponding basis functions. 

\par The aforementioned trend can also be observed in \reffig{fig:NNZ}, which illustrates the evolution of the number of non-zero coefficients (magnitude greater than $\delta_{rs}$ for RS) for both the LS and RS coefficients. The total number of coefficients corresponding to an $8^{th}$ order approximation turns out to be 1287. As time progresses, the LS tends to employ almost all these coefficients and calculates the smallest two-norm error, $\epsilon_{ls}$. The RS solution, on the other hand, only features around 50 coefficients by sacrificing accuracy for sparsity. 

\par To investigate the efficacy of the proposed scheme, the training and testing pdf errors are compared. The training data corresponds to the true pdf value at the CUT points, which are used to eventually compute the RS coefficients for the final approximation. The testing data corresponds to 10000 randomly generated MC points. \reffig{fig:TestTrainError} shows the averaged relative error between the true and approximated pdf values for both the training and the testing sets. The training error is low, as expected. The testing error is on the order of \num{e-2} after \SI{11}{\hour} of propagation, which roughly corresponds to about $1 \%$ of relative error at the most. This is a good approximation in the context of spacecraft tracking, as the errors incurred by the analytical approximation would still place the spacecraft in the expected region of search as would be obtained through numerical propagation of the true pdf. Therefore, the proposed approximation allows one to rapidly predict the state pdf of the spacecraft several hours in the future, by simply evaluating an analytical function facilitated by the coefficients computed a priori. 

\section{Summary}
    This paper is a follow up to a preceding paper that discussed the alternate CR3BP model, the S-VAM, in detail, and compared it with the traditional Cartesian model for its pros and cons \cite{sharan2022coordinate}. It was established that the S-VAM offers a way to quantify and propagate uncertainties in the direction of velocity independent of its magnitude, while strictly holding the Jacobi integral constant throughout the numerical propagation. Further study into the advantages of such a feature, and its exploitation in dealing with thrust uncertainties in the CR3BP is carried out in this work. As a result, several stochastic insights beneficial to spacecraft tracking are obtained from a completely different perspective compared to a traditional Cartesian view.
    
    \par To conclude, the development of the S-VAM in the preceding work \cite{sharan2022coordinate} laid a foundation for this paper to go forward and carry out detailed tests of the higher order STT computation scheme on a few different kinds of uncertainties associated with a thrusting maneuver. Furthermore, a strategy to obtain analytical approximations of the state pdfs at each time was also investigated successfully. 

    \par The analytical scheme developed in this work for pdf propagation is still in its infancy in the context of realistic trajectories in the CR3BP, and requires further fine-tuning to work for all trajectories for extended periods of time. Nevertheless, this work contains results that show great promise for this scheme, especially because it was tested on a scenario which had exaggerated values of uncertainty to begin with. Further development of this pdf propagation scheme would include process noise considerations and continuous thrust trajectories, which, in addition to the proven capability of the CUT-STTs in this work, would create a holistic framework to quantify and propagate uncertainties associated with cislunar trajectories with computational ease, while obtaining stochastic insights at high levels of accuracy. 

\section{Acknowledgement}
	This material is based upon work supported jointly by the AFOSR grants FA9550-20-1-0176 and FA9550-22-1-0092. 
    
\bibliographystyle{AAS_publication}   
\bibliography{references.bib}   

\end{document}